\documentclass[11pt]{amsart}

\usepackage{amsmath}
\usepackage{amssymb}
\usepackage{graphicx}

\hoffset = - 0.5in

\newtheorem{theorem}{Theorem}[section]

\newtheorem{proposition}[theorem]{Proposition}
\theoremstyle {definition}
\newtheorem{definition}[theorem]{Definition}

\newtheorem{remark}[theorem]{Remarks}

\numberwithin{equation}{section}

\renewcommand{\geq}{\geqslant}
\renewcommand{\leq}{\leqslant}

\title{More on Property($M$) in separable Banach spaces: Revised}

\author[Tim Dalby]{Tim Dalby}

\date{\today}

\keywords{fixed point property, property($M$), Kadec-Klee norm, weak normal structure}

\subjclass[2010]{46B10, 47H09, 47H10}

\email{tim\_dalby@bigpond.com}

\begin{document}

\parindent = 0pt
\parskip = 8pt

\begin{abstract}

Property($M$) in separable Banach spaces has played an important role in metric fixed point theory. This paper explores some of the Banach space properties that can be associated with Property($M$) and Property($M^*$).

\end{abstract}

\maketitle

\section{Introduction}

Property($M$) was introduced  by Kalton in [14].  This was in the context of M-ideals of compact operators. In that paper and a further one by Kalton and Werner, [15], various properties of Banach spaces that resulted from Property($M$) were presented.  Garcia-Falset and Sims, [10], took this idea and showed that Property($M$) implied the weak Fixed Point Property.  In [3], Dalby continued with this exploration and revealed how some properties, including Property($M$), were related.  This present paper looks again at what happens when a Banach space has Property($M$) or, in the dual, Property($M^*$). 

Because the Fixed Point Property (FPP) and weak Fixed Point Property (wFPP) are separably determined, all the Banach spaces are separable.

\section{Definitions}
 
An important Banach space property is Property($M$) where  whenever $x_n \rightharpoonup 0 $
 then $\limsup_{n \rightarrow \infty}\| x_n - x \|$ is a function of $\| x \|$ only.  An equivalent definition is the following.

\begin{definition}
Kalton [14]

A Banach space $X$ has Propery($M$) if whenever $x_n \rightharpoonup 0 \mbox{ and } \newline \| u \| \leq \| v \| \mbox{ then } \limsup_{n \rightarrow \infty} \| x_n + u \| \leq \limsup_{n \rightarrow \infty} \| x_n + v \|.$ 

\end{definition}

Property($M^*$) in $X^*$ is the same except that it involves weak* null sequences. 

A related property is WORTH. 
\newpage
\begin{definition}
Sims [23]

A Banach space, $X$, has WORTH if for every weak null sequence and every $x \in X$
 \[ \lim_{n \rightarrow \infty} | \| x_n - x \| - \| x_n + x \| |  = 0. \]

\end{definition}

This property was also introduced by Rosenthal [19] and, Cowell and Kalton [2].  The latter authors used the name Property(au).

Clearly, Property($M$) implies WORTH.

In $X^*$, WORTH* is the same but uses weak* null sequences.  In [17], Lima used the term Property($wM^*$) and in [2] it is Property(au*).  Dalby, [4], proved that WORTH* implies WORTH.  This result is also contained within the proof of lemma 2.2 of [20].

Another important Banach space property, in the context of fixed point theory, is Property($K$).  

\begin{definition}
Sims [24]

A Banach space $X$ has property($K$) if there exists $K \in [0, 1)$ such that whenever $x_n \rightharpoonup 0, \lim_{ n \rightarrow \infty} \| x_n \| = 1 \mbox{ and }  \liminf_{n \rightarrow \infty}\| x_n - x \|  \leq 1 \mbox{ then } \| x \| \leq K.$

\end{definition}

Note that implicit in this definition is the fact that $X$ cannot have the Schur propery.

Sims in [24] showed that Property($K$) implied weak normal structure and hence the wFPP.

In $X^*$ the property is called Property$(K^*$) and involves weak* null sequences.

A key property in metric fixed point theory and also in other areas of Banach space research is Opal's condition.

\begin {definition}

Opial [21] 

A Banach space has Opial's condition if
\[ x_n \rightharpoonup 0 \ \mbox {and } x \not = 0 \mbox { implies } \limsup_n \| x_n \| < \limsup_n \| x_n - x \|. \]
The condition remains the same if both the $\limsup$s are replaced by $\liminf$s.

If the $<$ is replaced by $\leq$ then the condition is called nonstrict Opial.  WORTH implies nonstrict Opial and hence so does Property($M$).
In $X^*$ the conditions are the same except that the sequence is weak* null.

\end{definition}

Later a modulus was introduced to gauge the strength of Opial's condition and a stronger version of the condition was defined.

\begin{definition}

Lin, Tan and Xu, [18]

Opial's modulus is
\[ r_X(c) = \inf \{ \liminf_{n \rightarrow \infty} \| x_n - x \| - 1: c \geq 0, \| x \| \geq c,  x_n \rightharpoonup 0 \mbox{ and }\liminf_{n \rightarrow \infty} \| x_n \| \geq 1 \}. \]
$X$ is said to have uniform Opial's condition if $r_X(c) > 0$ for all $c > 0.$  See [18] for more details.

Details of other properties such as weak normal structure and Kadec-Klee norm can be found in Kirk and Sims [16] or Goebel and Kirk [11].
\end{definition}

\section{Results}

The following theorem brings together, in a more coherent way,  various result that appear in [3].

\begin{theorem}
Let $X$ be a separable Banach space with Property($M$) then the following are equivalent.
\begin{enumerate}

\item[(a)] $c_0 \not \hookrightarrow X$

\item[(b)] $X$ has weak normal structure

\item[(c)] $X$ has Property($K$)

\item[(d)] $X$ has a Kadec-Klee norm

\item[(e)] $r_X(1) > 0.$

\end{enumerate}

\end{theorem}

\begin{proof}

(a) $\iff$ (b) This is corollary 3.5 of [3].

(a) $\iff$ (c)  This is a combination of theorem 3.3 of [3] and proposition 3.7 of [5].

(a) $\Rightarrow$ (d) This is theorem 3.2 of [3].

(c) $\iff$ (e) This is theorem 3.10 of [3].

(d) $\Rightarrow$ (c) Assume that $X$ has a Kadec-Klee norm but does not have Property($K$). Then there exists $x_n \rightharpoonup 0, \| x_n \| \rightarrow 1, (y_n) \mbox{ such that } \| y_n \| \rightarrow 1^{-} \mbox{ where }$
\[ \limsup_{n \rightarrow \infty} \| x_n - y_m \| \leq 1 \mbox{ for all } m. \]
Using the nonstrict Opial condition, which is implied by Property($M$), 
\[ 1 = \lim_{n \rightarrow \infty} \| x_n \|  = \limsup_{n \rightarrow \infty} \| x_n \| \leq \limsup_{n \rightarrow \infty} \| x_n - y_m \|  \leq 1 \mbox{ for all } m. \]

Thus $\limsup_{n \rightarrow \infty} \| x_n - y_m \|  = 1 \mbox{ for all } m.$

 Since $\| y_n \| \rightarrow 1^-$, Property($M$) and continuity of $\limsup_{n \rightarrow \infty} \| x_n - y \|$ together imply $\limsup_{n \rightarrow \infty} \| x_n - y \| = 1 \mbox{ for all } y \in B_X.$
 
 If $\| y \| = 1, \mbox{ then } x_n - y\rightharpoonup -y \mbox{ and } \limsup_{n \rightarrow \infty} \| x_n - y \|  = 1 = \| -y \|.$  The Kadec-Klee norm means $x_n - y\rightarrow -y \mbox{ and } x_n \rightarrow 0,$ this  contradicts $\| x_n \| \rightarrow 1.$
 
\end{proof}

A careful reading of the proofs presented in [3] and the one above  shows that the proofs readily transfer to $X^*$, Property($M^*$) and weak* null sequences.  For clarity and completeness the following theorem states what happens in $X^*.$

\begin{theorem}
Let $X$ be a separable Banach space with $X^*$ having Property($M^*$) then the following are equivalent.
\begin{enumerate}

\item[(a)] $c_0 \not \hookrightarrow X^*$

\item[(b)] $X^*$ has weak* normal structure

\item[(c)] $X^*$ has Property($K^*$)

\item[(d)] $X^*$ has a weak* Kadec-Klee norm

\item[(e)] $r_{X^*}(1) > 0.$

\end{enumerate}

\end{theorem}

In [14] Kalton showed that if $X^*$ has Property($M^*$) then $X$ has Property($M$).  The reverse implication requires a little more from $X$ and results in a little more in $X^*$.  Kalton and Werner, in [15], showed that if $X$ is a separable Banach space with  Property($M$) and $\ell_1 \not \hookrightarrow X$ then $X^*$ has Property($M^*$).  In addition, $X^*$ is separable and has a weak* Kadec-Klee norm. In light of the above theorem this means that $X^*$ satisfies the properties in that list.  Note dual spaces that are separable cannot contain $c_0$ so $c_0 \not \hookrightarrow X^*,$ which is a nice indication that the list is correct.  The following theorem summarises these points.

\begin{theorem}
Let $X$ be a separable Banach space with Property($M$) and  $\ell_1 \not \hookrightarrow X.$  Then $X^*$ has the following properties.
\begin{enumerate}

\item [(a)] Property($M^*$)

\item[(b)] $X^*$ is separable

\item[(c)] $c_0 \not \hookrightarrow X^*$

\item[(d)] $X^*$ has weak* normal structure

\item[(e)] $X^*$ has Property($K^*$)

\item[(f)] $X^*$ has a weak* Kadec-Klee norm

\item[(g)] $r_{X^*}(1) > 0.$

\end{enumerate}

\end{theorem}

There is a little bit more to this story.  In [17] Lima defined Property($wM^*$), which is the same as WORTH*, and showed that with this property,  $X^*$ has the Radon-Nikodym Property (RNP).  So every separable subspace of $X$ has a separable dual. Because, in this setting, $X$ is separable,  we have $X^*$ separable.  Thus $c_0 \not \hookrightarrow X^*$ even without the stronger Property($M^*$). This also means that  $\ell_1 \not \hookrightarrow X.$  For more background to this topic of the RNP see, for example, Diestel and Uhl [6]

This is summarised in the following proposition.

\begin{proposition}
Let $X$ be a separable Banach space where $X^*$ has WORTH* then $X^*$ has the following properties.

\begin{enumerate}

\item [(a)] Radon-Nikodym Property

\item[(b)] $X^*$ is separable

\item[(c)] $c_0 \not \hookrightarrow X^*$

\end{enumerate}

\end{proposition}

Also, $X^*$ having the RNP is equivalent $X$ being an Asplund space.  See, for example, van Dulst and Namioka [9].  So as stated above, $\ell_1 \not \hookrightarrow X.$  This leads to the following proposition.

\begin{proposition}
Let $X$ be a separable Banach space then $X^*$ has Property($M^*$) if and only if $X$ has Property($M$) and $\ell_1 \not \hookrightarrow X.$ 
\end{proposition}

It is not known if $X$ having WORTH and $\ell_1 \not \hookrightarrow X$ implies $X^*$ has WORTH*.  Cowell and Kalton in [2] tried but as they wrote "we do prove a result very close to it."   See theorem 6.4 of [2].

Combining theorem 3.2 and propositions 3.4 and 3.5 we have the following.

\begin{theorem}
Let $X$ be a separable Banach space where $X^*$ has Property($M^*$), or equivalently $X$ has Property($M$) and $\ell_1 \not \hookrightarrow X$, then $X^*$ has the following properties.

\begin{enumerate}

\item[(a)] $X^*$ is separable

\item[(b)] $c_0 \not \hookrightarrow X^*$

\item[(c)] $X^*$ has weak* normal structure

\item[(d)] $X^*$ has Property($K^*$)

\item[(f)] $X^*$ has a weak* Kadec-Klee norm

\item[(f)] $r_{X^*}(1) > 0.$

\end{enumerate}

\end{theorem}

In [15], lemma 3.5,  Kalton and Werner proved items (a) and (f).  In [1], in corollary 3.6, Cabello and Nieto proved item(c).

The presence or absence of either $c_0$ or $\ell_1$ play important roles in this area of Banach space theory.  Which brings to mind that in a Banach lattice, reflexivity is equivalent to not containing sublattices (or subspaces) isomorphic to $c_0$ and $\ell_1.$  A natural question is, is there  something similar happening if the space has Property($M$) or Property($M^*$)?  The following answers that question in the affirmative. This theorem formed part of a research report, written by the author,  from The University of New England, Australia in 2002 and a talk given at the Australian Mathematical Society annual conference at The University of Newcastle in 2002.

\begin{proposition}
Let $X$ be a separable Banach space.  If $X^*$ has Property($M^*$) then the following are equivalent.

\begin{enumerate}

\item [(a)] $X$ is reflexive

\item[(b)] $\ell_1 \not \hookrightarrow X^*$

\item[(c)] $c_0 \not \hookrightarrow X$

\item[(d)] $X$ has the FPP

\item[(e)] $X^*$ has the FPP.

\end{enumerate}

\end{proposition}

\begin{proof}
 $X$ being reflective is equivalent to $X^*$ being reflexive.  So (a) $\Rightarrow$ (b) and (a) $\Rightarrow$ (c) are well known.
 
 (b) $\Rightarrow$ (a)  In [14] Kalton proved that if $X^*$ has Property($M^*$) then $X$ is $M$-embedded.  That is, $X$ is a $M$-ideal in $X^{**}.$  So $X^*$ is $L$-embedded, see Harmand and Lima [12] or  Harmand, Werner and Werner [13].  Because $\ell_1 \not \hookrightarrow X^*, X^*$ must be reflexive.  See for example corollary 2.3 in Chapter IV of [13].
 
 (c) $\Rightarrow$ (a) Combining $X$ being $M$-embedded and $c_0 \not \hookrightarrow X,$ with corollary 3.7 in Chapter III of [13] proves $X$ is reflexive.
 
 (a) $\Rightarrow$ (d) Property($M^*$) implies $X$ has Property($M$).  Garcia-Falset and Sims, [10], have shown that Property($M$) leads to the wFPP.  Since $X$ is reflexive this is equivalent to the FPP.
 
(a) $\Rightarrow$ (e)  From proposition 3.6, $X^*$ has weak* normal structure. Theorem 4.3 of [11] states that in this case $X^*$ has the w*FPP.  Now reflexivity can be used to state that $X^*$ has the FPP.
 
 $\neg$ (a) $\Rightarrow \neg$ (d)  Assume that $X$ is nonreflexive then using proposition 6 from [22], because $X$ is $M$-embedded, it must contains an asymptotically isometric copy of $c_0$ and so $X$ does not have the FPP.  That last step relies on a result from [8].
 
 $\neg$ (a) $\Rightarrow \neg$ (e) Assume that $X$ is nonreflexive so $X^*$ is nonreflexive.  It is also $L$-embedded so using corollary 4 from [22], $X^*$ contains an asymptotically isometric copy of $\ell_1$ and so $X^*$ does not have the FPP. That last step relies on a result from [7].
\end{proof}

\begin{remark}

\begin{enumerate}

\item[(a)] If $X$ has Property($M$) and $\ell_1 \not \hookrightarrow X$ then the same list of equivalent conditions is true. 

\item[(b)] If  $X$ has Property($M$) then $X$ is reflexive if and only if $\ell_1 \not \hookrightarrow X$ and  $c_0 \not \hookrightarrow X.$  So, in this respect, $X$ is acting the same as a Banach lattice.

\item[(c)] A Banach space with Property($M^*$) joins a select collection of spaces that are reflexive if and only if they have the FPP.  Subspaces of $L_1[0,1]$ are in this class.

\item[(d)] If $X^*$ has Property($M^*$) then $X$  having a subspace isomorphic to $c_0$ is equivalent to $X$ contains an asymptotically isometric copy of $c_0$.  This is because if  $c_0  \hookrightarrow X$ then $X$ is nonreflexive and Pfitzner, [22], has shown that nonreflexive $M$-embedded Banach spaces  contain $c_0$ asymptotically isometrically.

\item[(e)] Similarly, if $X^*$ has Property($M^*$) then $X^*$ having a subspace isomorphic to $\ell_1$ is equivalent to $X^*$ containing $\ell_1$ asymptotically isometrically.  This comes from Pfitzner, [22], again, where since $X^*$ is nonreflexive and $L$-embedded it must contain $\ell_1$ asymptotically isometrically.
\end{enumerate}

\end{remark}


\begin{thebibliography}{99}

\bibitem {1} J. Cabello and E. Nieto, {\it On $M$-type structures and the fixed point property}, Houston J. Math. {\bf 26} (2000), 549-560.
\bibitem{2} S. Cowell and N. Kalton, {\it Asymptotic unconditionality}, Quarterly J. Math. {\bf 61} (2010), 217-240.
\bibitem {3} T. Dalby, {\it Relationships between properties that imply the weak fixed point property}, J. Math. Anal. Appl. {\bf 253} (2001), 578-589.
\bibitem {4} T. Dalby, {\it The effect of the dual on a Banach space and the weak fixed point property}, Bull. Austral. Math. Soc. {\bf 67} (2003), 177-185.
\bibitem{5} T. Dalby and B. Sims, {\it Banach lattices and the weak fixed point property}, Proceedings of the Seventh International Conference on Fixed Point Theory and its Applications, Guanjauato, Mexico (2005), 63-71.
\bibitem {6}J. Diestel and J. Uhl, Jr, {\it The Radon-Nikodym theorem for Banach space valued measures}, Rocky Mountain J. Math. {\bf 6} (1976), 1-46.
\bibitem {7} P. Dowling and C. Lennard, {\it Every nonreflexive subspace of $L_1[0,1]$ fails the fixed point property}, Proc. Amer. Math. Soc. {\bf 125} (1997), 443-446.
\bibitem {8} P. Dowling, C. Lennard and B. Turrent, {\sl Reflexivity and the fixed point property for nonexpansive maps}, J. Math. Anal. Appl. {\bf 200} (1996), 635-662.
\bibitem {9} Dulst, D. van and I. Namioka, {\sl A note on trees in conjugate Banach spaces.}, Indag. Math. (Proceedings) {\bf 87} (1984), 7-10.
\bibitem {10} J. Garc\'{i}a-Falset and B. Sims, {\it Property(M) and the weak fixed point property}, Proc. Amer. Math. Soc. {\bf 125} (1997), 2891-1896.
\bibitem {11} K. Goebel and W. Kirk, {\it Topics in metric fixed point theory}, Cambridge University Press, Cambridge, 1990.
\bibitem {12} P. Harmand and \r{A}.  Lima, {\it Banach spaces which are $M$-ideals in their bidual}, Trans. Amer. Math. Soc. {\bf 283} (1994), 253-264.
\bibitem {13} P. Harmand, D. Werner, and W. Werner, {\it M-ideals in Banach Spaces and Banach Algebras}, Lecture Notes in Math. {\bf 154} Springer,  Berlin, Heidelberg 1993.
\bibitem {14} N. Kalton, {\it M-ideals of compact operators}, Illinois J. Math. {\bf 208} (1993), 147-169.
\bibitem {15} N. Kalton and D. Werner, {\it Property(M), M-ideals, and almost isometric structure of Banach spaces}, J. Reine Angew. Math. {\bf 461} (1995), 137-178.
\bibitem {16} W. Kirk and B. Sims (ed.), {\it Handbook of metric fixed point theory}, Kluwer Academic Publishers, Dordrecht,  2001.
\bibitem {17} \r{A}. Lima, {\it Property(wM*) and the unconditional metric compact approximation property}, Studia Math. {\bf 113} (1995), 249-263.
\bibitem {18} P.-K. Lin, K.-K.Tan and H.-K. Xu, {\it Demiclosedness principle and asymptotic behavior for asymptotically nonexpansive mappings}, Nonlinear Anal. {\bf 24} (1995), 929-946.
\bibitem {19} H. Rosenthal, {\it Some remarks concerning unconditional basic sequences}, Longhorn Notes, Texas Functional Analysis Seminar 1982-1983 (1983), The University of Texas at Austin, 15-47.
\bibitem {20} S. Prus, {\it Banach spaces with the uniform Opial condition}, Nonlinear Anal. {\bf 18} (1992), 697-704.
\bibitem{21} Z. Opial, {\it Weak convergence of the sequence of successive approximations for nonexpansive mappings}, Bull. Amer. Math. Soc. {\bf 73} (1967), 591-597.
\bibitem {22} H. Pfitzner, {\it A note on asymptotically isometric copies of $\ell^1$ and $c_0$}, Proc. Amer. Math. Soc. {\bf 129} (2001), 1367-1373.
\bibitem {23} B. Sims, {\it Orthogonality and fixed points of nonexpansive maps}, Proc. Centre Math. Anal. Austral. Nat. Univ. {\bf 20} (1988), 178-186.
\bibitem {24} B. Sims, {\it A class of spaces with weak normal structure}, Bull. Austral. Math. Soc. {\bf 50} (1994), 523-528.

\end{thebibliography}
\end{document}